\documentclass[12pt]{amsart}

\usepackage{amssymb,amsthm,amsmath,hyperref,amsrefs,verbatim}

\usepackage{fullpage}
\newtheorem{theorem}{Theorem}
\newtheorem{corollary}{Corollary}
\newtheorem{lemma}{Lemma}
\newtheorem{proposition}{Proposition}
\newtheorem{definition}{Definition}

\theoremstyle{remark}
\newtheorem{remark}[theorem]{Remark}

\newcommand{\N}{\mathbb{N}}

\renewcommand{\epsilon}{\varepsilon}
\renewcommand{\leq}{\leqslant}

\title{Nuclear dimension and $n$-comparison}
\author{Leonel Robert}
\address{Department of Mathematics and Statistics, York University, Toronto, Canada M3J 1P3}
\email{lrobert@mathstat.yorku.ca}

\begin{document}

\begin{abstract}
It is shown that if a C$^*$-algebra has nuclear dimension $n$ then its Cuntz semigroup
has the property of $n$-comparison. It then follows from results by Ortega, Perera, and R\o rdam
that $\sigma$-unital C$^*$-algebras of finite nuclear dimension (and even of nuclear dimension at most $\omega$)
are stable if and only if they have no non-zero unital quotients and no non-zero bounded traces. 
\end{abstract}

\maketitle

\section{Introduction}

In \cite{winter-zacharias1}, Winter and Zacharias define nuclear dimension for C$^*$-algebras.
This is a form of noncommutative dimension which directly generalizes the covering dimension of topological spaces.
Finite nuclear dimension is specially relevant to the classification of C$^*$-algebras. The simple C$^*$-algebras of finite nuclear dimension have been proposed as a likely class of C$^*$-algebras for which Elliott's classification in terms of K-theory and traces holds true.

In the main result of this paper it is shown that the Cuntz semigroup of a C$^*$-algebra of finite nuclear dimension $n$ satisfies the $n$-comparison property. For $n=0$, this property is the same as almost unperforation in the Cuntz semigroup. For arbitrary $n$,
it is reminiscent of the comparability  between vector bundles whose fibrewise dimensions differ sufficiently relative to the dimension of the base space.

The $n$-comparison property for the Cuntz semigroup
was first considered by Toms and Winter (see \cite{toms-winter}*{Lemma 6.1}). They showed that $n$-comparison holds under the more restrictive assumption that the
C$^*$-algebra  is simple unital of  decomposition rank $n$ (the decomposition rank bounds the nuclear dimension, and, unlike nuclear dimension, 
is infinite for UCT Kirchberg algebras).
The $n$-comparison property was subsequently studied, and more precisely defined,  by Ortega, Perera and R\o rdam in \cite{ortega-perera-rordam}. These authors obtained a  simple criterion of stability for $\sigma$-unital C$^*$-algebras with  $n$-comparison in their Cuntz semigroup: the  C$^*$-algebra is stable if and only if it has no non-zero unital quotients and no non-zero bounded 2-quasitraces. By Theorem \ref{main} below, this stability criterion
applies to all C$^*$-algebras of finite nuclear dimension. It then follows that $\sigma$-unital C$^*$-algebras of finite nuclear dimension have the corona factorization property.

Let us recall the definition of nuclear dimension given in \cite{winter-zacharias1}.

\begin{definition}\label{nc}
The C$^*$-algebra $A$ has nuclear dimension $n$ if $n$ is the smallest natural number
for which there exist nets of completely positive contractions (henceforth abbreviated as c.p.c.) 
\[
\psi_\lambda^i\colon A\to F_\lambda^i\hbox{ and  }\phi_\lambda^i\colon F_\lambda^i\to A,
 \]
with  $i=0,1,\dots,n$, $\lambda\in \Lambda$, and $F_\lambda^i$ finite dimensional C$^*$-algebras for all $i$ and $\lambda$, such that

(i) $\phi_\lambda^i$ is an order 0 map (i.e., preserves orthogonality) for all $i$ and $\lambda$,

(ii) $\lim_\lambda\sum_{i=1}^n \phi_\lambda^i\psi_\lambda^i(a)=a$ for all $a\in A$.

If no such $n$ exists then $A$ has infinite nuclear dimension.
\end{definition}

Let us recall the definition given in \cite{ortega-perera-rordam} of the $n$-comparison property of an ordered semigroup.
For $x,y$ elements of an ordered semigroup $S$, let us write $x\leq_s y$ if $(k+1)x\leq ky$ for some $k\in \mathbb{N}$.

\begin{definition}
The ordered semigroup $S$ has the $n$-comparison property if  $x\leq_s y_i$ for $x,y_i\in S$ and $i=0,1,\dots,n$, implies $x\leq \sum_{i=0}^n y_i$.
\end{definition}

Let $Cu(A)$ denote the stabilized Cuntz semigroup of the C$^*$-algebra $A$ (i.e., the semigroup  $W(A\otimes \mathcal K)$; see \cite{coward-elliott-ivanescu}).

It is shown in Lemma \ref{tracial} below that for $Cu(A)$ the $n$-comparison property  can be reformulated 
as follows: if $[a],[b_i]\in Cu(A)$, with $i=0,1,\dots,n$, satisfy that 
for each $i$ there is $\epsilon_i>0$ such that $d_\tau(a)\leq (1-\epsilon_i)d_\tau(b_i)$ for all dimension functions $d_\tau$ induced by lower
semicontinuous 2-quasitraces,
then $[a]\leq \sum_{i=0}^n [b_i]$. It is this formulation of the $n$-comparison property 
that is used by Toms and Winter in \cite{toms-winter}, and that may potentially have the most applications.

\begin{theorem}\label{main}
If $A$ has nuclear dimension  $n$ then $Cu(A)$ has the $n$-comparison property.
\end{theorem}

The following section is dedicated to the proof of Theorem \ref{main}. The last section discusses 
the application of Theorem \ref{main}, and of a variation on it that relates to  $\omega$-comparison, to 
establishing the  stability of C$^*$-algebras of finite (or at most $\omega$) nuclear dimension.

\section{Proof of Theorem \ref{main}}
Let us start by proving that the property of $n$-comparison may be formulated using comparison 
by lower semicontinuous 2-quasitraces instead of the relation $\leq_s$. This result, however, will not be needed
in the proof of Theorem \ref{main}. 

For $[a],[b]\in Cu(A)$, elements of the Cuntz semigroup of $A$,
let us write $[a]<_\tau[b]$ if there is $\epsilon>0$ such that $d_\tau(a)\leq (1-\epsilon)d_\tau(b)$ for all dimension functions $d_\tau$ induced by lower semicontinuous 2-quasitraces $\tau\colon A^+\to [0,\infty]$
(see \cite[Section 4]{elliott-robert-santiago}). We do not assume that the 2-quasitraces are necessarily finite on a dense subset of $A^+$.

\begin{lemma}\label{tracial}
The ordered semigroup $Cu(A)$ has the property of $n$-comparison if and only if
for $[a],[b_i]\in Cu(A)$, with $i=0,1,\dots,n$, $[a]<_\tau [b_i]$ for all $i$ implies that 
$[a]\leq \sum_{i=0}^n [b_i]$. 
\end{lemma}

\begin{proof}
It is clear that if $[a]\leq_s [b]$ then $[a]<_\tau [b]$. Thus, $n$-comparison is implied by the property
stated in the lemma. Suppose that we have $n$-comparison in $Cu(A)$. Let $[a],[b_i]\in Cu(A)$, with $i=0,1,\dots,n$, be such that  $[a]<_\tau [b_i]$ for all $i$. Let us show that $[(a-\epsilon)_+]\leq_s [b_i]$ for all $\epsilon>0$ and all $i$. Let $\lambda\colon Cu(A)\to [0,\infty]$
be additive and order preserving and let us define $\tilde\lambda\colon Cu(A)\to [0,\infty]$ by 
$\tilde \lambda([c])=\sup_{\delta>0} \lambda([(c-\delta)_+])$. It is known that there is a lower semicontinuous
2-quasitrace $\tau$ such that $\tilde \lambda([a])=d_\tau(a)$ for all $a\in (A\otimes \mathcal K)^+$ (see \cite[Proposition 4.2, Lemma 4.7]{elliott-robert-santiago}). 
We have 
\[
\lambda([(a-\epsilon)_+])\leq \tilde\lambda([a])\leq (1-\epsilon_i)\tilde \lambda([b_i])\leq 
(1-\epsilon_i)\lambda([b_i])
\]
for all $i$. By \cite[Proposition 2.1]{ortega-perera-rordam}, this implies that $[(a-\epsilon)_+]\leq_s [b_i]$ for all $i$. Since $Cu(A)$ has $n$-comparison, $[(a-\epsilon)_+]\leq \sum_{i=0}^n [b_i]$. Taking supremum over $\epsilon>0$, we get that $[a]\leq \sum_{i=0}^n [b_i]$.
\end{proof}

Throughout the rest of this section $\Lambda$ denotes the index set in the definition of nuclear dimension. This set may be chosen
to be the pairs $(F,\epsilon)$ where $F\subseteq A$ is finite and $\epsilon>0$. Let us denote by $A_\Lambda$ the algebra $\prod_{\lambda} A\slash \bigoplus_\lambda A$. Let $\iota\colon A\to A_\Lambda$ denote the diagonal embedding of $A$ into $A_\Lambda$.

\emph{Notation convention.} For a family of C$^*$-algebras $(A_\lambda)$ we use the notation $\bigoplus_\lambda A_\lambda$ to refer
to the C$^*$-algebra of nets $(x_\lambda)$ such that $\|x_\lambda\|\to 0$, while $\prod_\lambda A_\lambda$ denotes the nets of uniformly bounded norm. 

In \cite[Proposition 3.2]{winter-zacharias1} Winter and Zacharias  show 
that if $A$ has nuclear dimension $n$ then the maps $\psi_\lambda^i$ in the definition of nuclear dimension 
may be chosen asymptotically of order 0. That is, such that the induced maps $\psi^i\colon A\to \prod_\lambda F_\lambda^i/\bigoplus_\lambda F_\lambda^i$, for $i=0,1,\dots,n$, have order 0 .
We get the following proposition as a result of this.

\begin{proposition}\label{ncrestated}
If $A$ has nuclear dimension $n$ then for $i=0,1,\dots,n$ there are c.p.c. order 0 maps 
$\psi^i\colon A\to \prod_\lambda F_\lambda^i/\bigoplus_\lambda F_\lambda^i$ and $\phi^i\colon \prod_\lambda F_\lambda^i/\bigoplus _\lambda F_\lambda^i\to A_\Lambda$
such that
\begin{align}\label{iota}
\iota=\sum_{i=0}^n \phi^i\psi^i,
\end{align}
where $\iota\colon A\to A_\Lambda$ is the diagonal embedding of $A$ into $A_\Lambda$.
\end{proposition}
\begin{proof}
As pointed out in the previous paragraph, by \cite[Proposition 3.2]{winter-zacharias1} the maps $\psi_\lambda^i$ in the definition
of nuclear dimension may be chosen so that the induced maps $\psi^i$ are of order 0.

The equation \eqref{iota} is a  consequence of Definition \ref{nc} (ii).

Let us show that the maps  $\phi^i\colon \prod_\lambda F_\lambda^i/\bigoplus_\lambda F_\lambda^i\to A_\Lambda$, induced by the order 0 maps $\phi_\lambda^i$,
also have order 0. 
It is clear that $(\phi^i_\lambda)\colon \prod_\lambda F_\lambda^i\to \prod_\lambda A$ has order 0. Hence, 
$\alpha\circ (\phi^i_\lambda)\colon \prod_\lambda F_\lambda^i\to A_\Lambda$, where $\alpha$ is the quotient
onto $A_\Lambda$, has order 0. We will be done once we show that if $\phi\colon C\to D$ is a c.p.c. map of  order 0
and $\phi|_I=0$ for some closed two sided ideal $I$, then the induced map $\tilde \phi\colon C/I\to D$
has order 0. By \cite[Corollary 3.1]{winter-zacharias2}, there is a *-homomorphisms $\pi\colon C\otimes C_0(0,1]\to D$ such that $\phi(c)=\pi(c\otimes t)$ for all $c\in C$.
From $\pi(I\otimes t)=0$ we get that  $\pi(I\otimes C_0(0,1])=0$. Thus, $\pi$ induces a *-homomorphism $\tilde \pi\colon C/I\otimes C_0(0,1]\to D$.
Since $\tilde \phi(c)=\tilde \pi(c\otimes t)$ for all $c\in C/I$,  $\tilde \phi$ has order 0.
\end{proof}

An ordered semigroup $S$ is said unperforated if $kx\leq ky$
for $x,y\in S$ and $k\in \N$ implies $x\leq y$.

\begin{lemma}\label{al-un}
(i) If $Cu(A)$ is unperforated then so is $Cu(A/I)$ for any closed two-sided ideal $I$.

(ii) If $(A_\lambda)_\lambda$ are C$^*$-algebras such that  $Cu(A_\lambda)$ is unperforated for all $\lambda$
then so are $Cu(\prod_\lambda A_\lambda)$ and $Cu(\prod_\lambda A_\lambda/\bigoplus_\lambda A_\lambda)$.
\end{lemma}

\begin{proof}
(i) Let $[\tilde a],[\tilde b]\in Cu(A/I)$ be such that $k[\tilde a]\leq k[\tilde b]$ for some $k\in \mathbb{N}$.
Then for $[a],[b]\in Cu(A)$ lifts of $[\tilde a]$ and $[\tilde b]$ we have 
\[
k[a]\leq k[b]+[i]\leq k([b]+[i]), 
\]
for some $i\in (I\otimes \mathcal K)^+$ (by \cite{ciuperca-robert-santiago}*{Proposition 1}). Since $Cu(A)$ is  unperforated, we have $[a]\leq [b]+[i]$, and passing
to $Cu(A/I)$ we get that $[\tilde a]\leq [\tilde b]$.

(ii) Let $(a_\lambda)_\lambda,(b_\lambda)_\lambda\in (\prod_\lambda A_\lambda)\otimes \mathcal K\subseteq \prod_\lambda (A_\lambda \otimes \mathcal K)$
be positive elements  of norm at most 1 such that $k[(a_\lambda)_\lambda]\leq k[(b_\lambda)_\lambda]$ for some $k$. Let $\epsilon>0$. Then there is $\delta>0$
such that 
\[
k[((a_\lambda-\epsilon)_+)_\lambda]\leq k[((b_\lambda-\delta)_+)_\lambda]. 
\]
Hence, $k[(a_\lambda-\epsilon)_+]\leq k[(b_\lambda-\delta)_+]$ for all $\lambda$. Since $Cu(A_\lambda)$ is  unperforated,
 $[(a_\lambda-\epsilon)_+]\leq [(b_\lambda-\delta)_+]$. 
Let $x_\lambda\in A_\lambda\otimes \mathcal K$ be such that 
\[
(a_\lambda-2\epsilon)_+=x_\lambda^*x_\lambda\hbox{ and  }x_\lambda x_\lambda^*\in \mathrm{her}((b_\lambda-\delta)_+). 
\]
Since $\|x_\lambda\|^2\leq \|a_\lambda\|\leq 1$, we have 
$(x_\lambda)_\lambda\in \prod_\lambda (A_\lambda\otimes \mathcal K)$. We also
have $x_\lambda x_\lambda^*\leq \frac{1}{\delta}b_\lambda$, whence 
\[
((a_\lambda-2\epsilon)_+)_{\lambda}=(x_\lambda)_\lambda^*(x_\lambda)_\lambda\hbox{ and }
(x_\lambda)_\lambda(x_\lambda)_\lambda^*\in \mathrm{her}((b_\lambda)_\lambda). 
\]
But  $(\prod_\lambda A_\lambda)\otimes \mathcal K$ sits as a hereditary subalgebra of $\prod_\lambda (A_\lambda \otimes \mathcal K)$. Therefore, $(x_\lambda)_\lambda\in (\prod_\lambda A_\lambda)\otimes \mathcal K$ and so $[((a_\lambda-2\epsilon)_+)_\lambda]\leq [(b_\lambda)_\lambda]$. Since $\epsilon$ may be arbitrarily small, we have
$[(a_\lambda)_\lambda]\leq [(b_\lambda)_\lambda]$ as desired.
\end{proof}

\begin{lemma}\label{iotalem}
Let $a,b\in (A\otimes \mathcal K)^+$. If $[\iota(a)]\leq [\iota(b)]$ then $[a]\leq [b]$.
\end{lemma}

\begin{proof}
Let  $\epsilon>0$. Since $[\iota(a)]\leq [\iota(b)]$,  there is $d\in A_\Lambda\otimes \mathcal K$
such that $d^*\iota(b)d=\iota(a-\epsilon)_+$. That is,
there are $d_\lambda\in A_\lambda\otimes \mathcal K$ such that $d_\lambda^* b d_\lambda\to (a-\epsilon)_+$. Thus,
$[(a-\epsilon)_+]\leq [b]$. Since $\epsilon>0$ may be  arbitrarily small, we have $[a]\leq [b]$.   
\end{proof}

\begin{remark}
In proving Theorem \ref{main}, a stronger property than $n$-comparison
will be shown to hold for C$^*$-algebras of nuclear dimension $n$: if $x,y_i\in Cu(A)$, with $i=0,1,\dots,n$,
satisfy that $k_ix\leq k_iy_i$ for some $k_i\in \mathbb{N}$ and all $i$, then $x\leq \sum_{i=0}^n y_i$. 
This property, unlike $n$-comparison, does not seem to have a formulation in terms
of comparison by lower semicontinuous 2-quasitraces. 
\end{remark}

\begin{proof}[Proof of Theorem \ref{main}]
Suppose that there are $k_i\in \mathbb{N}$ such that $k_i[a]\leq k_i[b_i]$ for $i=0,1,\dots,n$.
Since c.p.c. order 0 maps preserve Cuntz comparison (by \cite[Corollary 3.5]{winter-zacharias2}), we have that 
$k_i[\psi^i(a)]\leq k_i[\psi^i(b_i)]$ for all $i$. Since the Cuntz semigroup of finite dimensional
algebras is unperforated, we have by Lemma \ref{al-un} that 
the Cuntz semigroup of $\prod_\lambda F_\lambda^i/\bigoplus_\lambda F_\lambda^i$ is unperforated. Thus, 
$[\psi^i(a)]\leq [\psi^i(b_i)]$. The maps $\phi^i$ preserve Cuntz equivalence (since they are c.p.c. of order 0), 
whence
\[
[\phi^i\psi^i(a)]\leq [\phi^i\psi^i(b_i)]\leq [\sum_{j=0}^n \phi^j\psi^j(b_i)]=[\iota(b_i)].
\]
So,
\[
[\iota(a)]=[\sum_{i=0}^n \phi^i\psi^i(a)]\leq \sum_{i=0}^n [\phi^i\psi^i(a)]\leq \sum_{i=0}^n [\iota(b_i)].
\]
By Lemma \ref{iotalem}, this implies that $[a]\leq \sum_{i=0}^n [b_i]$.
\end{proof}

\section{Stability of C$^*$-algebras}
A stable C$^*$-algebra has no non-zero unital quotients and no non-zero bounded 2-quasitraces 
(see \cite{ortega-perera-rordam}*{Proposition 4.6}). 
In \cite[Proposition 4.8]{ortega-perera-rordam} Ortega, Perera and R\o rdam show that the converse is true provided that 
the C$^*$-algebra is $\sigma$-unital and its Cuntz semigroup has the $n$-comparison property. This, combined with Theorem \ref{main}
and the fact that for exact C$^*$-algebras bounded 2-quasitraces are traces, implies that a $\sigma$-unital C$^*$-algebra of finite nuclear dimension is
stable if an only if it has no non-zero unital quotients and no nonzero bounded traces. 
Ortega, Perera, and R\o rdam also show  that $\omega$-comparison, a weakening of $n$-comparison, suffices to obtain the same
stability criterion.

\begin{definition}(c.f. \cite{ortega-perera-rordam}*{Definition 2.9})
Let $S$ be an ordered semigroup closed under the suprema of increasing sequences. Then $S$ has the $\omega$-comparison property if $x\leq_s y_i$, for $x,y_i\in S$ and $i=0,1,\dots$, implies $x\leq \sum_{i=0}^\infty y_i$.
\end{definition} 

\begin{remark} The definition of $\omega$-comparison given above differs slightly from the definition given in \cite{ortega-perera-rordam}. Nevertheless, both definitions agree for ordered semigroups in the category $Cu$ introduced in \cite{coward-elliott-ivanescu}, and therefore, also for ordered semigroups arising as Cuntz semigroups of C$^*$-algebras. 
\end{remark}

A notion of nuclear dimension at most $\omega$ may be modelled after the statement of Proposition \ref{ncrestated}.

\begin{definition}\label{atmostomega}
Let us say that a C$^*$-algebra $A$ has nuclear dimension at most $\omega$ if for $i=0,1,2\dots$ there are nets of c.p.c maps 
$\psi_\lambda^i\colon A\to F_\lambda^i$ and
$\phi_\lambda^i\colon F_\lambda^i\to A$,  with $F_\lambda^i$ finite dimensional C$^*$-algebras, such that

(ii) the induced maps $\psi^i\colon A\to \prod_\lambda F_\lambda^i/\bigoplus_\lambda F_\lambda^i$ and 
$\phi^i\colon \prod_\lambda F_\lambda^i/\bigoplus_\lambda F_\lambda^i\to A_\Lambda$ are c.p.c. order 0,

(iii) $\iota(a)=\sum_{i=0}^\infty \phi^i\psi^i(a)$ for all $a\in A$ (the sum on the right side is norm convergent).  
\end{definition}
For example, if the C$^*$-algebras $(A_i)_{i=0}^\infty$ all have finite nuclear dimension, then
$\bigoplus_{i=0}^\infty A_i$ has nuclear dimension at most $\omega$. It is not clear whether the 
assumption that the maps $\psi_\lambda^i$ be asymptotically order 0 may be dropped in Definition \ref{atmostomega} (and then proved),  or  if the other results on finite nuclear dimension proved in \cite{winter-zacharias1}
also hold for nuclear dimension at most $\omega$. 

The proof of Theorem \ref{main} goes through,  mutatis mutandis, for nuclear dimension at most $\omega$.
We thus have,
\begin{theorem}\label{main2}
If $A$ has nuclear dimension at most $\omega$ then $Cu(A)$ has the $\omega$-comparison property.
\end{theorem}
 
Combined with the results of \cite{ortega-perera-rordam}, Theorem \ref{main2} yields the following corollary, which improves
 on \cite{ng-winter}*{Theorem 0.1}, \cite{ortega-perera-rordam}*{Corollary 4.10} and \cite{ortega-perera-rordam}*{Corollary 5.13}.
\begin{corollary}
Let $A$ be a C$^*$-algebra of  nuclear dimension at most $\omega$ and let $B\subseteq A\otimes\mathcal K$ be hereditary and $\sigma$-unital. 
Then $B$ is stable if and only if it has no non-zero unital quotients and no non-zero bounded traces. 
$B$ has the corona factorization  property.
 
\end{corollary}
\begin{proof}
By Theorem \ref{main2}, $Cu(A)$ has the $\omega$-comparison property. Since $Cu(B)$ is an ordered semigroup  ideal of $Cu(A)$,
$Cu(B)$ has the $\omega$-comparison property too. Hence, by \cite{ortega-perera-rordam}*{Proposition 4.8}, $B$ is stable if and only if
it has no non-zero unital quotients and no non-zero bounded traces. 

Let $C\subseteq B\otimes \mathcal K$ be full and hereditary, and suppose that $M_n(C)$ is stable. Then $M_n(C)$, and consequently $C$, cannot have non-zero unital quotients or bounded traces. Thus $C$ is stable. This shows that $B$ has the corona factorization property.
\end{proof}

\begin{remark} In the hypotheses of \cite{ortega-perera-rordam}*{Proposition 4.8} the C$^*$-algebra
$A$ is assumed separable. A closer look into the proof of this result reveals that it suffices to assume
that  the hereditary subalgebra $B\subseteq A\otimes \mathcal K$ be $\sigma$-unital. This justifies the application
of \cite{ortega-perera-rordam}*{Proposition 4.8} in the above proof.
\end{remark}

\textbf{Acknowledgment.} I am grateful to Wilhelm Winter for pointing out a correction to the proof of Lemma \ref{al-un}.

\end{document}